\documentclass[11pt]{amsart}

\usepackage{latexsym}
\usepackage{amssymb}
\usepackage{amscd}
\usepackage[all]{xy}
\usepackage[mathscr]{euscript}
\usepackage{dsfont}

\normalsize

\RequirePackage{xspace}

\newcommand{\thought}[1]{}
\renewcommand{\thought}[1]{ \textbf{[#1]}}

\usepackage{enumerate}        
\newenvironment{roenumerate}{\begin{enumerate}[\upshape (i)]}{\end{enumerate}}

\newcommand\nc {\newcommand}
\newcommand\rnc{\renewcommand}

\newcount\blopone
\newcount\xone
\newcount\xtwo
\newcount\ytwo

\newtheorem{theorem}{Theorem}[section]
\newtheorem{prop}[theorem]{Proposition}
\newtheorem{conjecture}[theorem]{Conjecture}
\newtheorem{com}[theorem]{Comment}
\newtheorem{apl}[theorem]{Application}
\newtheorem{exercise}[theorem]{Exercise}
\newtheorem{redu}[theorem]{Reduction}
\newtheorem{refinement}[theorem]{Refinement}
\newtheorem{summary}[theorem]{Summary}
\newtheorem{importnota}[theorem]{Important Notation}
\newtheorem{prblm}[theorem]{Problem}
\newtheorem{notation}[theorem]{Notation}
\newtheorem{explanation}[theorem]{Explanation}
\newtheorem{defin}[theorem]{Definition}
\newtheorem{caution}[theorem]{Caution}
\newtheorem{remark}[theorem]{Remark}
\newtheorem{reminder}[theorem]{Reminder}
\newtheorem{illustration}[theorem]{Illustration}
\newtheorem{observation}[theorem]{Observation}
\newtheorem{lemma}[theorem]{Lemma}
\newtheorem{construction}[theorem]{Construction}
\newtheorem{discussion}[theorem]{Discussion}
\newtheorem{corollary}[theorem]{Corollary}
\newtheorem{example}[theorem]{Example}
\newtheorem{conclusion}[theorem]{Conclusion}
\newtheorem{sketch}[theorem]{Sketch}
\newtheorem{triviality}[theorem]{Triviality}
\newtheorem{proto}[theorem]{Prototype Quasifibration}
\newtheorem{cauex}[theorem]{Cautionary Example}
\newtheorem{hypo}[theorem]{Hypothesis}
\newtheorem{subth}{ }[theorem]
\newtheorem{case}{Case}[theorem]
\newtheorem{ssubth}{ }[subth]
\newtheorem{facts}[theorem]{Facts}
\newtheorem{history}[theorem]{Historical Survey}
\newtheorem{proofs}[theorem]{Discussion of the Proofs, Old and New}

\nc\tri[1]{\begin{triviality}
\label{#1}}
\nc\cnj[1]{\begin{conjecture}
\label{#1}}
\nc\fac[1]{\begin{facts}
\label{#1}
\begin{em}}
\nc\app[1]{\begin{apl}
\label{#1}
\begin{em}}
\nc\skt[1]{\begin{sketch}
\label{#1}
\begin{em}}
\nc\hst[1]{\begin{history}
\label{#1}
\begin{em}}
\nc\pfs[1]{\begin{proofs}
\label{#1}
\begin{em}}
\nc\cas[1]{\begin{case}
\label{#1}
\begin{em}}
\nc\rfn[1]{\begin{refinement}
\label{#1}}
\nc\prt[1]{\begin{proto}
\label{#1}}
\nc\lem[1]{\begin{lemma}
\label{#1}}
\nc\pro[1]{\begin{prop}
\label{#1}}
\nc\thm[1]{\begin{theorem}
\label{#1}}
\nc\dis[1]{\begin{discussion}
\label{#1}
\begin{em}}
\nc\cor[1]{\begin{corollary}
\label{#1}}
\nc\dfn[1]{\begin{defin}
\label{#1}}
\nc\sthm[1]{\begin{subth}
\label{#1}
\begin{em}}
\nc\exm[1]{\begin{example}
\label{#1}
\begin{em}}
\nc\obs[1]{\begin{observation}
\label{#1}
\begin{em}}
\nc\plm[1]{\begin{prblm}
\label{#1}
\begin{em}}
\nc\rmk[1]{\begin{remark}
\label{#1}
\begin{em}}
\nc\rmd[1]{\begin{reminder}
\label{#1}
\begin{em}}
\nc\ntn[1]{\begin{notation}
\label{#1}
\begin{em}}
\nc\exe[1]{\begin{exercise}
\label{#1}
\begin{em}}
\nc\xpl[1]{\begin{explanation}
\label{#1}
\begin{em}}
\nc\smr[1]{\begin{summary}
\label{#1}
\begin{em}}
\nc\cau[1]{\begin{caution}
\label{#1}
\begin{em}}
\nc\hyp[1]{\begin{hypo}
\label{#1}
\begin{em}}
\nc\imn[1]{\begin{importnota}
\label{#1}
\begin{em}}
\nc\rdn[1]{\begin{redu}
\label{#1}
\begin{em}}
\nc\cax[1]{\begin{cauex}
\label{#1}
\begin{em}}
\nc\cmt[1]{\begin{com}
\label{#1}
\begin{em}}
\nc\con[1]{\begin{construction}
\label{#1}
\begin{em}}
\nc\ill[1]{\begin{illustration}
\label{#1}
\begin{em}}
\nc\ssthm[1]{\begin{ssubth}
\label{#1}
\begin{em}}
\nc\cnc[1]{\begin{conclusion}
\label{#1}
\begin{em}}

\nc\elem{\end{lemma}}
\nc\erdn{\end{em}\end{redu}}
\nc\erfn{\end{refinement}}
\nc\eprt{\end{proto}}
\nc\ethm{\end{theorem}}
\nc\ecor{\end{corollary}}
\nc\edfn{\end{defin}}
\nc\esthm{
\end{em}
\end{subth}}
\nc\epro{\end{prop}}
\nc\etri{\end{triviality}}
\nc\ecnj{\end{conjecture}}
\nc\eexm{\end{em}
\end{example}}
\nc\eobs{\end{em}
\end{observation}}
\nc\ecmt{\end{em}
\end{com}}
\nc\efac{\end{em}
\end{facts}}
\nc\eapp{\end{em}
\end{apl}}
\nc\ermk{\end{em}
\end{remark}}
\nc\ermd{\end{em}
\end{reminder}}
\nc\eill{\end{em}
\end{illustration}}
\nc\eplm{\end{em}
\end{prblm}}
\nc\ecas{\end{em}
\end{case}}
\nc\eskt{\end{em}
\end{sketch}}
\nc\ecau{\end{em}
\end{caution}}
\nc\ecax{\end{em}
\end{cauex}}
\nc\eimn{\end{em}
\end{importnota}}
\nc\entn{\end{em}
\end{notation}}
\nc\eexe{\end{em}
\end{exercise}}
\nc\expl{\end{em}
\end{explanation}}
\nc\edis{\end{em}
\end{discussion}}
\nc\econ{\end{em}
\end{construction}}
\nc\esmr{\end{em}
\end{summary}}
\nc\ehst{\end{em}
\end{history}}
\nc\epfs{\end{em}
\end{proofs}}
\nc\ehyp{
\end{em}
\end{hypo}}
\nc\ecnc{\end{em}
\end{conclusion}}
\nc\essthm{\end{em}
\end{ssubth}}

\nc\sst{\scriptstyle}
\newcommand{\comment}[1]{}
\newcommand{\ri}{\longrightarrow}

\newcommand{\zz}{{\mathbb Z}}

\newcommand{\K}{{\mathbf K}}

\newcommand{\D}{{\mathbf D}}

\newcommand{\rr}{{\mathbb R}}
\newcommand{\C}{{\mathbb C}}
\newcommand{\pp}{{\mathbb P}}

\nc\op{^{\hbox{\rm\tiny op}}}
\nc\mth{^{\hbox{\rm\tiny th}}}

\nc\script{\mathscr}
\nc\z{\zeta}
\nc\bc{{\mathbb{BC}}}
\nc\KK{{\mathbb{K}}}
\nc\ct{{\script T}}
\nc\cf{{\script F}}
\nc\cg{{\script G}}
\nc\ch{{\script H}}
\nc\ck{{\script K}}
\nc\cl{{\script L}}
\nc\cm{{\script M}}
\nc\cn{{\script N}}
\nc\cv{{\script V}}
\nc\ce{{\script E}}
\nc\cs{{\script S}}
\nc\car{{\script R}}
\nc\cd{{\script D}}
\nc\cc{{\script C}}
\nc\ca{{\script A}}
\nc\ci{{\script I}}
\nc\cj{{\script J}}
\nc\co{{\script O}}
\nc\cu{{\script U}}
\nc\cw{{\script W}}
\nc\cx{{\script X}}
\nc\Cp{{\script P}}
\nc\cq{{\script Q}}
\nc\cy{{\script Y}}
\nc\cz{{\script Z}}
\nc\bd{\begin{description}}
\nc\ed{\end{description}}
\nc\ctob{{\script C}at\big(\ci^{op},\ca\big)}
\nc\clim{{\ds\mathop{\rm lim}_{\ds\longleftarrow}}\,}
\nc\climi{\clim_{\!i}\,}
\nc\climn{\clim^{\!n}\,}
\nc\colim{{\ds\mathop{\rm colim}_{\ds\la}}}
\nc\colimj{{\ds\mathop{\rm colim}_{\ds\la}}{}_{j\,}}
\nc\oa{\overline{\ca}}
\nc\s{\sigma}
\nc\ta{\tau}
\nc\os{\overline\sigma}
\nc\ot{\overline\tau}
\nc\T{\Sigma}
\nc\Tm{\Sigma^{-1}}
\nc\de[1]{{\mathop{\rm deg(#1)}}}
\nc\Ad[1]{\mathop{\rm Ad}(#1)}
\nc\ad[1]{\mathop{\rm ad}(#1)}
\nc\kth{{\it K}--theory}
\nc\loc[1]{{\text{\rm Loc}(#1)}}
\nc\coloc[1]{{\text{\rm Coloc}(#1)}}

\def\der #1 {D\left(#1\right)}
\nc\prf{\begin{proof}}
\nc\eprf{\end{proof}}
\nc\ds{\displaystyle}
\nc\Tor{\text{\rm Tor}}

\nc\cb{{\script B}}
\nc\ab{{\script A}b}

\nc\be{\begin{roenumerate}}
\nc\ee{\end{roenumerate}}

\nc\cat[1]{{\script C}at\Big({\big\{#1\big\}}\op\,\,,\,\,\ab\Big)}
\nc\csab{{\script C}at\big(\cs^{op},\ab\big)}
\nc\ctab{{\script C}at\Big({\{\ct^\alpha\}}^{op},\ab\Big)}
\nc\csex{{\script E}x\big(\cs^{op},\ab\big)}
\nc\ctex{{\script E}x\Big({\{\ct^\alpha\}}^{op},\ab\Big)}
\nc\sub{\qquad\subset\qquad}
\nc\ctr[1]{{\left.\ct\left(-,#1\right)\right|}_{\cs}}
\nc\ctrf[2]{{\left.\ct\left(#1,#2\right)\right|}_{\cs}}
\nc\Ctr[1]{{\left.\ct\left(-,#1\right)\right|}_{\ct^\alpha}}
\nc\Ctrf[2]{{\left.\ct\left(#1,#2\right)\right|}_{\ct^\alpha}}

\nc\la{\longrightarrow}
\nc\nin{\noindent}
\nc\cad[1]{\text{card}(#1)}
\nc\eq{\quad=\quad}
\nc\BA{\begin{array}{c}}
\nc\EA{\end{array}}
\nc\barr{
\[
\begin{array}{cccccccccccccccc}
}
\nc\earr{
\end{array}
\]
}
\nc\as[1]{{\langle S\rangle}^{#1}}
\nc\sh{\text{\it shift}}

\nc\yy[1]{{\left.\ct\left(-,#1\right)\right|}_{\ct^c}}
\nc\vrep[2]{{\left.\ct\left(#1,#2\right)\right|}_{\ct^\alpha}}
\nc\da{\downarrow}
\nc\Hom{{\mathop{\rm Hom}}}
\nc\HHom{{\script H}{\mathop{\rm om}}}
\nc\End{{\mathop{\rm End}}}
\nc\Ext{{\mathop{\rm Ext}}}
\nc\mr\Modtc
\nc\PExt{{\mathop{\rm PExt}}}
\nc\stm{\text{\rm stmod}(kG)}
\nc\stM{\text{\rm StMod}(kG)}
\nc\e{\varepsilon}
\nc\p{\varphi}

\nc\rs{\s^{-1}A}
\nc\br{{\{\s^{-1}A\}}}
\nc\ra\ri
\nc\y[1]{\mathbf{y}#1}
\nc\x[1]{\mathbf{z}#1}
\nc\mmod[1]{\text{\rm mod--}#1}
\nc\Mod[1]{#1\text{--\rm Mod}}
\nc\Md {\ensuremath{\mathop{\textup{Mod}}}}
\rnc\mod[1]{\ensuremath{\mathop{#1\textup{--mod}}}\xspace}
\nc\MMod[1]{\text{Mod-}#1}
\nc\Modtc{\Mod{\ct^c}}
\nc\pgldim[1]{\mathop{\rm pgldim}\,#1}
\nc\tf{{\rm [TR5]}}
\nc\tfs{{\rm [TR5$^*$]}}
\nc\Fun{\text{\rm Funct}(F\op,\ab)}
\nc\sym{\text{\rm Sym}}
\nc\sgn{\text{\rm sgn}}
\nc\Pro{\text{\rm Prod}^{}_\alpha(F\op,\ab)}
\nc\Yt[1]{{\left.\Hom_\ct^{}\left(-,#1\right)\right|}_F^{}}
\nc\dl{\delta}
\nc\Proj[1]{#1\text{--\rm Proj}}
\nc\proj[1]{#1\text{--\rm proj}}
\nc\Flat[1]{#1\text{--\rm Flat}}
\nc\Inj[1]{#1\text{--\rm Inj}}
\nc\Ima{\mathrm{Im}}
\nc\Ker{\mathrm{Ker}}
\nc\ov{\overline}
\nc\wt{\widetilde}
\nc\wh{\widehat}
\nc\ph{\varphi}
\nc\tstr{{\it t}--structure}
\nc\tstrs{{\it t}--structures}
\nc\spec[1]{{\text{\rm Spec}(#1)}}

\nc\EProd{\text{\rm EProd}}
\nc\ECoprod{\text{\rm ECoprod}}
\nc\Prod{\text{\rm Prod}}
\nc\ldimp{\text{\rm LDim}^{\prod}}
\nc\ldimc{\text{\rm LDim}^{\coprod}}
\nc\gen[2]{{\langle#1\rangle}^{}_{#2}}
\nc\Gen[2]{{\big\langle#1\big\rangle}^{}_{#2}}
\nc\genu[3]{{\langle#1\rangle}^{[#3]}_{#2}}
\nc\ogen[1]{\ov{\langle#1\rangle}}
\nc\ogenun[2]{\ov{\langle#1\rangle}_{#2}^{}}
\nc\ogenu[3]{\ov{\langle#1\rangle}^{[#3]}_{#2}}
\nc\ogenul[3]{\ov{\langle#1\rangle}^{(-\infty,#3]}_{#2}}
\nc\ogenuf[3]{\ov{\langle#1\rangle}^{[#3,\infty)}_{#2}}
\nc\genuf[3]{{\langle#1\rangle}^{[#3,\infty)}_{#2}}
\nc\genul[3]{{\langle#1\rangle}^{(-\infty,#3]}_{#2}}
\nc\dperf[1]{\D^{\mathrm{perf}}(#1)}
\nc\dperfs[2]{\D_{#1}^{\mathrm{perf}}(#2)}
\nc\dcoh{\mathbf{D}^b_{\mathrm{coh}}}
\newcommand{\Dqc}{{\mathbf D_{\text{\bf qc}}}}
\newcommand{\Dqcs}[1]{{\mathbf D_{\text{\bf qc},#1}}}

\nc\dmcoh{\mathbf{D}^-_{\mathrm{coh}}}
\nc\dscoh{\mathbf{D}^{}_{\mathrm{coh}}}
\nc\RHHom{{\script{RH}}{\mathrm{om}}}
\nc\Coprod{\mathrm{Coprod}}
\nc\COprod{\mathrm{coprod}}
\nc\add{\mathrm{add}}
\nc\Add{\mathrm{Add}}
\nc\Smr{\mathrm{smd}}
\nc\id{\mathrm{id}}
\nc\LL{\mathbf{L}}
\nc\R{\mathbf{R}}
\nc\CC{\mathbf{C}}
\nc\wi{\wt{\text{\it\i}}}
\nc\exal{\ce\text{\it x}_\alpha(\ct^\alpha,\ab)}
\nc\exalz{\ce\text{\it x}_{\aleph_0}^{}(\ct^\alpha,\ab)}

\nc\tst[1]{\left({#1}^{\leq0},{#1}^{\geq0}\right)}
\nc\perf[1]{\script{P}\mathrm{erf}\left(#1\right)}
\nc\SF[1]{\mathrm{SF}\left(#1\right)}

\nc\one{\mathds{1}}
\nc\fc{\mathfrak{C}}
\nc\fl{\mathfrak{L}}
\nc\fs{\mathfrak{S}}
\nc\Prf{\text{\bf Perf}}
\nc\qc{\text{\bf qc}}
\nc\coh[1]{\script{C}\mathit{oh}\left(#1\right)}
\nc\vect[1]{\script{V}\mathit{ect}\left(#1\right)}
\nc\vectgr[1]{\script{V}\mathit{ect}\script{G}\mathit{r}\left(#1\right)}
\nc\Ch{\mathbf{Ch}}
\nc\fgt{\mathbf{fgt}}

\nc\hoco{
\begin{picture}(40,10)
\put(20,0){\makebox(0,0)[b]{\text{\rm Hocolim}}}
\put(5,-2){\vector(1,0){30}}
\end{picture}\,\,}

\nc\holim{{
\begin{picture}(40,10)
\put(20,0){\makebox(0,0)[b]{\text{\rm Holim}}}
\put(35,-2){\vector(-1,0){30}}
\end{picture}}}
\renewcommand\geq{\geqslant}
\renewcommand\leq{\leqslant}

\begin{document}

\author{Amnon Neeman}\thanks{The research was partly supported 
by the Australian Research Council}
\address{Centre for Mathematics and its Applications \\
        Mathematical Sciences Institute\\
        Building 145\\
        The Australian National University\\
        Canberra, ACT 2601\\
        AUSTRALIA}
\email{Amnon.Neeman@anu.edu.au}

\title{Obstructions to the existence of Bounded {\it t}--structures}

\begin{abstract}
In a striking 2019 article, Antieau, Gepner and Heller found 
{\it K--}theoretic obstructions
to bounded t-structures.
 We will survey
their work, as well as some progress since. The focus will be
on the open problems that arise from this.
\end{abstract}

\subjclass[2010]{Primary 18G80, secondary 19D35, 14F08}

\keywords{Derived categories,
bounded {\it t}--structures, \kth, schemes, vector bundles}

\maketitle

\tableofcontents

\setcounter{section}{-1}

\section{Introduction}
\label{S0}

The concept of \tstr s on
triangulated categories, which is
recalled in Definition~\ref{D3.1},
was introduced by
Be{\u\i}linson, Bernstein and 
Deligne~\cite[Chapter~1]{BeiBerDel82}, in their attempt
to better understand MacPherson's work on intersection
homology. In the following decades \tstr s rapidly gained 
attention, and were 
successfully adapted and employed as a 
technique applicable across much of 
mathematics. And bounded t-structures
have been of particular interest. 
\exm{E0.1}
To list a few examples of subjects where bounded
\tstr s have played a pivotal role:
\sthm{E0.1.1}
In the study of the representations of finite groups of Lie
type, Deligne and Lusztig~\cite{Deligne-Lusztig79} 
introduced the powerful theory of character sheaves.
It has come to play an enormous role in the field, and
at its center (at least in modern developments of
the theory) are bounded
perverse \tstr s, on derived categories of sheaves 
over certain 
varieties.
\esthm
\sthm{E0.1.2}
The theory of stability conditions on triangulated
categories, whose origins go back to 
Bridgeland~\cite{Bridgeland07,Bridgeland08}, has turned
out to be of great interest
in algebraic geometry. It is also growing in prominence
in the field of geometric representation theory. And one
way to view a 
stability condition is as a family of 
bounded \tstr s, indexed
by the set $\rr$ of real numbers, and satisfying a 
list of properties we omit.
\esthm
\sthm{E0.1.3}
Motivic homotopy theory, whose origins go back to a string of
deeply influential articles by Voevodsky, was developed 
as an attempt to move forward on Grothendieck's vision
of a motivic cohomology theory. More precisely it was
inspired by the Be{\u\i}linson-Bloch approach to the
problem. Anyway: in the decades that followed
Grothendieck's formulation of his program, no one had made
any substantial progress on constructing the 
conjectural abelian
category $\ca$ of motives. Voevodsky's idea was that $\D^b(\ca)$,
its bounded derived category, should be constructible by
homotopy-theoretic methods. And the search for $\ca$ was 
transformed into the search for a bounded \tstr\ on $\D^b(\ca)$,
with heart $\ca$.
\esthm
\eexm
The conventional wisdom, acquired over three and a half 
deacades, was that \tstr s are plentiful and diverse. There is 
a plethora of substantially different  
techniques for constructing them.
\exm{E0.3}
Below is a partial list of known methods to produce bounded \tstr s:
\sthm{E0.3.1}
If a triangulated category $\cs$ is built up of two triangulated categories
$\car$ and $\ct$ by recollement, then any pair of bounded \tstr s on $\car$
and on $\ct$ glue to a bounded \tstr\ on $\cs$. This technique
goes back to the birth of the subject in
Be{\u\i}linson, Bernstein and 
Deligne~\cite[Chapter~1]{BeiBerDel82},
where it was used to create the perverse \tstr s on certain
derived categories
of sheaves.
\esthm
\sthm{E0.3.2}
One can tilt a \tstr\ with respect to a torsion pair, a
technique introduced
in the work of Happel, Reiten and Smal\o~\cite{Happel-Reiten-Smalo96}.
This has been extensively studied since.
\esthm
\sthm{E0.3.3}
Given a set $S$ of objects in a triangulated category $\ct$, one
can look for the minimal \tstr\ that the set $S$ generates. By this
we mean: the \tstr\ $\big(\ct^{\leq0},\ct^{\geq0}\big)$ where $S$ is 
contained in $\ct^{\leq0}$ and $\ct^{\leq0}$ is minimal subject to 
the condition. The idea was introduced by
Alonso, Jerem{\'{\i}}as and Souto~\cite{Alonso-Jeremias-Souto03},
and extended by several authors since.

These \tstr s are rarely bounded, but can restrict to bounded \tstr s
on triangulated subcategories of $\ct$. A careful study of such examples can be
found, for example, in Alonso, Jerem{\'{\i}}as and 
Saor{\'{\i}}n~\cite{Alonso-Jeremias-Saorin10}.
\esthm
\sthm{E0.3.4}
One can construct bounded \tstr s on $\ct$ 
using Harder-Narasimhan filtrations.
This process begins with a central charge, which is a 
group homomorphism $Z:K_0(\ct)\la\C$. And then, given any object $X\in\ct$,
one inductively finds the most destabilizing subobject of
$X$ to produce a filtration. This technique leads to a Bridgeland 
stability condition. The original paper introducing
the method was Bridgeland~\cite{Bridgeland07}, although
the reader might prefer the more elegant,
later version due to Bayer~\cite{Bayer19}. Anyway: at this point
the literature on the active subject is immense.
\esthm
\eexm

One can therefore imagine the surprise that greeted the
2019 article by Antieau, Gepner and
Heller~\cite{Antieau-Gepner-Heller19}, which finds 
{\it K--}theoretic obstructions to the existence of bounded
\tstr s. 
Since the aim of this survey is to highlight open questions,
let us insert one already.

\plm{P0.1}
In \ref{E0.1.3} above we mentioned that the search for a
good abelian category $\ca$
of motives, as envisioned by Grothendieck,
can be reformulated as the search for the triangulated
category $\D^b(\ca)$ and a bounded \tstr\ on it.
There are by now several candidates for what should be the right
triangulated category $\D^b(\ca)$, and in this active field
more candidates are being produced by the day. (OK, this might
be an exaggeration, but new candidates pop up with remarkable
frequency). 

In the light of Antieau, Gepner and
Heller~\cite{Antieau-Gepner-Heller19} we now know,
as mentioned above, that 
there are {\it K--}theoretic 
obstructions to the existence of bounded \tstr s.
The challenge to the experts 
becomes to compute these obstructions
for the many candidate $\D^b(\ca)$'s, and weed
the ones with no chance of working. This means: if an 
obstruction
is nonzero for a candidate $\D^b(\ca)$,
then this candidate can safely be discarded. Either there is
no way the putative $\D^b(\ca)$ 
could have a bounded \tstr, or else
any bounded \tstr\ on it has to have a heart which is
absurd enough to be ruled out as a potential
abelian category of motives.
\eplm

In this survey we will remind the reader of the 
{\it K--}theoretic background, explain in more
detail what 
Antieau, Gepner and
Heller~\cite{Antieau-Gepner-Heller19} proved and
conjectured, explain why some of these conjectures seemed
highly plausible while others represented a daring leap of
faith, describe the progress made on the 
conjectures since \cite{Antieau-Gepner-Heller19} appeared,
and (most importantly) highlight the many fascinating 
questions that remain.

\bigskip

\nin
{\bf Acknowledgements.}\ \ The author is grateful to Ben Antieau,
Jeremiah Heller and Evgeny Shinder for
improvements on earlier versions.

\section{A quick reminder of algebraic $K$-theory}
\label{S1}

In this section we give a brief and minimalistic overview of
algebraic \kth. 

In an abelian category it makes sense to speak of exact sequences. Recall:
given a category $\ca$, being abelian is a property, not a
structure. The category $\ca$ is abelian if it has finite limits and
colimits, and these satisfy some properties 
we do not repeat here. Thus if
a category $\ca$ happens to be abelian then the exact sequences in it
make sense intrinsically.

Now if $\ca$ is an abelian category and $\ce\subset\ca$ is a 
full subcategory, then it also makes sense to speak of exact sequences
in $\ce$. These can be defined to be those diagrams in $\ce$ which become
exact when viewed in $\ca$. Of course: this definition is no longer
intrinsic in $\ce$, it depends on our chosen embedding into the abelian
category $\ca$. The axiomatic characterization of those $\ce$'s
which admit embeddings into abelian categories,
as extension-closed full subcategories, is due
to Quillen~\cite{Quillen1}. We recall

\rmd{R1.1}
An \emph{exact category} is an additive category $\ce$, together with 
a prescribed class of \emph{admissible exact sequences}
\[\xymatrix{
A\ar[r]^-f & B\ar[r]^-g & C\ ,
}\]
satisfying a list of properties we omit. And it is a theorem that
a category $\ce$ is exact if and only if there exists a 
fully faithful functor
$F:\ce\la\ca$, where $\ca$ is an abelian category,
such that the essential image
of $F$ is extension-closed, and a sequence $A\la B\la C$ is admissible
in $\ce$ if and only if $0\la F(A)\la F(B)\la F(C)\la 0$ is exact
in $\ca$.
\ermd

\rmk{R1.3}
Note that in the exact category $\ce$ the admissible 
exact sequences are a
\emph{structure,} not a \emph{property.} They need to be specified in 
advance. An additive category $\ce$ can have more than one exact structure,
as we will see below.

Let $\cf$ be any additive category. We can turn $\cf$ into an exact
category by stipulating that the admissible exact sequences are the
``split exact sequences'', meaning the isomorphs of diagrams
\[\xymatrix{
A\ar[r]^-i & A\oplus C\ar[r]^-\pi & C
}\]
where $i$ is the inclusion and $\pi$ is the projection. We will denote
the additive category $\cf$, with the split exact structure, by the 
symbol $\cf^\oplus$. Given any exact category $\ce$, we can produce 
a possibly different exact structure by forming
$\ce^\oplus$. And we stress that normally 
one would not expect these two exact
structures to agree. This will play a key role in Section~\ref{S5},
starting with Discussion~\ref{D5.7}.
\ermk

\rmk{R1.5}
In Reminder~\ref{R1.1} we told the reader that we omit the
full list of
properties that the admissible exact sequences are required to satisfy.
But there
is one property worth mentioning, even in a minimalistic survey.

In an abelian category $\ca$, the diagram
\[\xymatrix{
0\ar[r] &A\ar[r]^-f & B\ar[r]^-g & C\ar[r] & 0
}\]
is a short exact sequence if and only if $f$ is the kernel of $g$ and
$g$ is the cokernel of $f$. In an exact category this is no longer
an ``if and only if'' statement, after all Remark~\ref{R1.3}
teaches us that there cannot exist a characterization of the
admissible short exact sequences purely in terms of categorical 
data from $\ce$. But it is true that, if the
diagram
\[\xymatrix{
A\ar[r]^-f & B\ar[r]^-g & C
}\]
is an admissible short exact sequence, then $f$ is the kernel of
$g$ and $g$ is the cokernel of $f$. The reader can see that this follows
immediately from the existence of an exact embedding $F:\ce\la\ca$ as
in Reminder~\ref{R1.1}.
\ermk

\rmd{R1.7}
Let $\ce$ be an
essentially small exact category, and let
$F(\ce)$ be the free abelian group on the isomorphism
classes $[A]$ of objects $A\in\ce$.
The \emph{Grothendieck group $K_0(\ce)$}
is defined by the formula
\[
K_0(\ce)\eq\frac{F(\ce)}{\langle[A]-[B]+[C]\text{ for every admissible }A\la B\la C\rangle}
\]
In words: $K_0(\ce)$ is the quotient
of $F(\ce)$ by the relations generated by expressions $[A]-[B]+[C]$, where
$A\la B\la C$ runs over the admissible exact sequences in $\ce$.
\ermd

\dfn{D1.9}
Let $X$ be a scheme. The special case where $\ce$ is the exact category
of vector bundles on $X$ will interest us enough to deserve a symbol:
we denote this exact category by
$\vect X$. The exact structure is that a sequence $\cv'\la \cv\la\cv''$,
of vector bundles on $X$, is declared to be admissible
exact if and only if it is exact
in the abelian category of all sheaves.
And we will adopt the abbreviation
\[
K_0(X):=K_0\big(\vect X\big).
\]
\edfn 

\rmd{R1.11}
For any morphism of schemes $f:X\la Y$, there is an induced pullback
map $f^*:\vect Y\la\vect X$. Since this is an exact functor (meaning it
respects admissible short exact sequences), it induces a group homomorphism
$K_0(Y)\la K_0(X)$.

Now suppose we are given a scheme $X$, and two open subsets $U,V\subset X$.
Then there is a commutative square of inclusions
\[\xymatrix{
  U\cap V \ar[r] \ar[d] & U \ar[d] \\
  V\ar[r] & X\ ,
}\]
which the functor $K_0$ takes to a commutative square of group
homomorphisms
\[\xymatrix{
  K_0(U\cap V)  & K_0(U) \ar[l] \\
  K_0(V)\ar[u] & K_0(X)\ar[l]\ar[u]\ .
}\]
We can ``fold'' this commutative square into 
\[\xymatrix{
K_0(X)\ar[r]^-f & K_0(U)\oplus K_0(V)\ar[r]^-g & K_0(U\cap V)
}
\]
where the composite $gf$ vanishes. And it is classical that, \emph{as
  long as $U$ and $V$ cover $X$} (meaning $X=U\cup V$), this sequence is
exact in the middle.

It becomes natural to try to extend to a long
exact sequence. And after much work by many people this was done.
Assuming $X$ satisfies the resolution property (meaning
that there are enough vector bundles on $X$),
there is a way to define positive and negative 
algebraic \kth\ in such a way
that the above extends to a long exact Mayer-Vietoris sequence
\[
\xymatrix{
 & & \ar[dll]\\
K_{1}(X) \ar[r]& K_{1}(U)\oplus K_{1}(V)\ar[r] & K_{1}(U\cap V)\ar[dll] \\
K_{0}(X) \ar[r]& K_{0}(U)\oplus K_{0}(V)\ar[r] & K_{0}(U\cap V)\ar[dll] \\
K_{-1}(X) \ar[r]& K_{-1}(U)\oplus K_{-1}(V)\ar[r] & K_{-1}(U\cap V)\ar[dll] \\
   &
}\]
\ermd

\rmk{R1.13}
It goes without saying that, in addition to having
long exact sequences, this higher
algebraic \kth\ satisfies
many other
good properties, which we omit. Its definition and
functoriality properties extend
well beyond vector bundles on
schemes, and even well beyond exact categories.

In this survey, we will use only a tiny portion of the
vast general theory. For us, the important observations
are:
\sthm{R1.13.1}
  It is possible to define this higher 
algebraic \kth, both positive
  and negative, for any category of cochain complexes with
  a \emph{model structure}\footnote{By
    \emph{model structure} we mean either
    a biWaldhausen complicial category in the sense of
    Thomason and Trobaugh~\cite{ThomTro}, or a
    stable infinity category
    as in Lurie~\cite{Lurie09}.}.
\esthm
\sthm{R1.13.2}
  The Mayer-Vietoris long exact sequence, of
  Reminder~\ref{R1.11}, is a special case of much more
  general long exact sequences which will come up later.
\esthm
\sthm{R1.13.3}
  To each model category as in \ref{R1.13.1} one can associate a
  homotopy category, which is a triangulated
  category. More formally: there is a functor $\fgt$,
  from the category of model categories to the category
  of triangulated categories.
  And there is a theorem saying that, if
  $F:\cm\la\cn$ is a morphism
  of model categories such that 
  $\fgt(F):\fgt(\cm)\la\fgt(\cn)$ is an equivalence
  of triangulated categories,
  then the maps
  $K_i(F):K_i(\cm)\la K_i(\cn)$ are isomorphisms
  for all $i\in\zz$.
\esthm
It is standard to refer to the model category $\cm$ as
an \emph{enhancement} of the triangulated category
$\fgt(\cm)$.
\ermk

\section{Vanishing theorems for negative \kth}
\label{S2}

Next we want to quickly remind the reader that, in
a wide range of 
important special cases, there are vanishing theorems
for negative \kth. And as in Section~\ref{S1} we begin
with the special case of schemes.

\rmd{R2.1}
In Section~\ref{S1} we left out the history: there
was no mention of how and when higher algebraic
\kth\ was
developed, and who contributed what to this achievement.
Without changing this substantially, let us
note that Chuck Weibel was one of the key players
in the study of negative \kth, with substantial
contributions spanning decades. And already in
the 1980 article~\cite[Question~2.9]{Weibel80} he
formulated what came to be known as Weibel's conjecture.
\begin{itemize}
\item
  If $X$ is a noetherian scheme of dimension $n<\infty$,
  then $K_i(X)=0$ for all $i<-n$.
\end{itemize}
In the following decades much work was done, proving
special cases of this conjecture. In full generality
the proof may be found in the 2018 article by
Kerz, Strunk and
Tamme~\cite[Theorem~B]{Kerz-Strunk-Tamme18}.
It took a long time to settle this question.
\ermd

\plm{P2.3}
It would be interesting to discover the right
generalization of this theorem. Presumably there should
be a theorem saying that, for some large class
of model categories $\cm$, we have $K_i(\cm)=0$ for
all $i\ll0$. And, like Weibel's conjecture,
maybe this general theorem will come with an effective
bound, giving an explicit integer $n$, depending on $\cm$, 
such that
$K_i(\cm)=0$ for
all $i<-n$.

I do have a candidate conjecture, but it
would require too much notation to introduce it here.
And, in any case, the reader is encouraged to use her
imagination and come up with conjectures of her own.
\eplm

In the remainder of the article we will focus on
a different set of vanishing conjectures for negative
\kth.

\rmd{R2.5}
By a classical result due to Bass it was known that
\be
\item
For regular,
noetherian, finite-dimensional, \emph{affine} 
schemes $X$ one has $K_i(X)=0$ for all $i<0$.
\setcounter{enumiv}{\value{enumi}}
\ee
Now
recall that, for any scheme and by Definition~\ref{D1.9},
$K_i(X)=K_i\big(\vect X\big)$. The category $\vect X$,
of vector bundles on the noetherian scheme $X$, admits
an exact embedding into the
abelian category $\coh X$ of
all coherent sheaves on $X$. And it is a classical
theorem due to Quillen
that
\be
\setcounter{enumi}{\value{enumiv}}
\item
For finite-dimensional, regular, noetherian
schemes $X$ with enough vector bundles, the natural maps
\[\xymatrix{
K_i\big(\vect X\big)\ar[r] &
K_i\big(\coh X\big)  
}\]
are isomorphisms for all $i\in\zz$.
\ee
\ermd

\rmd{R2.7}
Schlichting generalized the known
statements of Reminder~\ref{R2.5}. The general vanishing
statements he proved, for negative \kth, are
as follows:
\be
\item
  Let $\ca$ be an essentially small abelian category.
  Then $K_{-1}(\ca)=0$.
\item
  Let $\ca$ be a \emph{noetherian} abelian category.
  Then $K_i(\ca)=0$ for all $i<0$.
\setcounter{enumiv}{\value{enumi}}  
\ee
The reader can find (i) in
\cite[Theorem~6 of Section~10]{Schlichting06}, while (ii)
is \cite[Theorem~7 of Section~10]{Schlichting06}.
By applying (ii) above to the noetherian
abelian category $\coh X$, and combining with
Quillen's isomorphism of Reminder~\ref{R2.5}(ii), we obtain
\be
\setcounter{enumi}{\value{enumiv}}
\item
For finite-dimensional, regular, noetherian
schemes $X$ with enough vector bundles, 
the abelian groups $K_i(X)=
K_i\big(\vect X\big)$ vanish for all $i<0$.
\setcounter{enumiv}{\value{enumi}}
\ee
In other words: as an easy  corollary of his
results, Schlichting was able to
generalize the vanishing statement of
Reminder~\ref{R2.5}(i) from affine schemes to all
schemes. And then he went on to conjecture
that even more should be true.
\be
\setcounter{enumi}{\value{enumiv}}
\item
For any abelian category $\ca$ and all integers $i<0$,
the groups $K_i(\ca)$ vanish.
\ee
This may be found in
\cite[Conjecture~1 of Section~10]{Schlichting06}.
The statement (iv) came to be known as
Schlichting's conjecture.
\ermd

\rmk{R2.9}
Let us explain why Schlichting's conjecture is plausible.
There is an old result of
Quillen~\cite[Section~5, Theorem~4]{Quillen1} asserting:
\sthm{R2.9.1}
  Let $\cb$ be an abelian category, and assume
  $\ca\subset\cb$ is a Serre subcategory. Recall: this
  means that $\ca$ is an abelian subcategory of $\cb$
  such that every $\cb$-subquotient of an object $A\in\ca$
  belongs to $\ca$. Now write $\cc$ for the
  abelian quotient category $\cc=\cb/\ca$.

  Then there is a long exact sequence in \kth\ 
\[
\xymatrix@R-10pt@C+20pt{
 & & \ar[dll]\\
K_{2}(\ca) \ar[r]& K_{2}(\cb)\ar[r] & K_{2}(\cc)\ar[dll]  &\\
K_{1}(\ca) \ar[r]& K_{1}(\cb)\ar[r] & K_{1}(\cc)\ar[dll]  & \\
K_{0}(\ca) \ar[r]& K_{0}(\cb)\ar[r] & K_{0}(\cc)\ar[r] &0\\
}\]
\esthm

\nin
So the challenge becomes to try to show that this
long exact sequence is sufficient to define negative
\kth\ without ever having to leave the world of
abelian categories.

The standard methods of homological algebra tell us that
this will work as long as there are enough abelian
categories $\cb$ with vanishing \kth. If we could
embed every abelian category $\ca$ as a Serre subcategory
$\ca\subset\cb$, where the higher \kth\ of $\cb$ must
vanish for some formal reason, then we'd be in business.

Let me not elaborate here
on the technical difficulties that
people encountered in their attempts to do this.
\ermk

\section{A reminder of \tstr s}
\label{S3}

Before we recall the formal definitions, let us give some intuition. A
\tstr\ on a triangulated category amounts to specifying which objects
in the category we declare to be positive, and which objects in the category
we declare to be negative. And then some conditions must be satisfied.
More formally, we have

\dfn{D3.1}
Let $\ct$ be a triangulated category. A \tstr\ on $\ct$ is a pair
of full subcategories $\big(\ct^{\leq0},\ct^{\geq0}\big)$, with
$\ct^{\leq0}$ to be thought of as the \emph{negative objects} and
$\ct^{\geq0}$ to be thought of as the \emph{positive objects.} And these
must satisfy the conditions
\be
\item
$\T\ct^{\leq0}\subset\ct^{\leq0}$ and $\ct^{\geq0}\subset\T\ct^{\geq0}$.
\item
$\Hom\big(\T\ct^{\leq0},\ct^{\geq0}\big)=0$.
\item
For every object $B\in\ct$ there exists a triangle $A\la B\la C\la$ in 
$\ct$, with $A\in\T\ct^{\leq0}$ and $C\in\ct^{\geq0}$.
\setcounter{enumiv}{\value{enumi}}
\ee
So much for the general definition of \tstr. Now we furthermore define
\be
\setcounter{enumi}{\value{enumiv}}
\item
The \tstr\ $\big(\ct^{\leq0},\ct^{\geq0}\big)$ is declared
\emph{bounded} if, for every object $X\in\ct$, there exists an
integer $n>0$ with $\T^nX\in\ct^{\leq0}$ and with $\T^{-n}X\in\ct^{\geq0}$.
\setcounter{enumiv}{\value{enumi}}
\ee 
\edfn

The most classical example is

\exm{E3.3}
Let $R$ be an associative ring, and let $\D(R)$ be the derived category
whose objects are cochain complexes of left $R$-modules. Then the following
gives a \tstr\ on $\D(R)$:
\begin{eqnarray*}
\D(R)^{\leq0}&=& \{X\in\D(R)\mid H^i(X)=0\text{ for all }i>0\}\ , \\
\D(R)^{\geq0}&=& \{X\in\D(R)\mid H^i(X)=0\text{ for all }i<0\}\ . \\
\end{eqnarray*}
This \tstr\ is not bounded. If we let $X$ be the cochain complex
\[\xymatrix{
\cdots\ar[r] & R\ar[r]^-0 & R\ar[r]^-0 & R\ar[r]^-0 &
R\ar[r]^-0 & R\ar[r] & \cdots
}\]
where in every degree we put the rank-1 free $R$-module, and all
the differentials are zero, then $X$ is an object of $\D(R)$ but
no shift of it is either positive or negative.

But, when restricted to the full subcategory $\D^b(R)$, the the 
\tstr\ above becomes bounded. The category $\D^b(R)$ has for objects
the cochain complexes that vanish outside a bounded interval. 
\eexm

In Example~\ref{E0.3} we told the reader that there are many other
examples of bounded \tstr s, and listed four known, general procedures
that produce them.

There is one general theorem about \tstr s that we will frequently appeal
to, hence we remind the reader:

\thm{T3.5}
Let $\ct$ be a triangulated category, and assume that 
$\big(\ct^{\leq0},\ct^{\geq0}\big)$ is a \tstr\ on $\ct$.
Then the category $\ct^\heartsuit=\ct^{\leq0}\cap\ct^{\geq0}$ is 
abelian. It is called the \emph{heart} of the \tstr.
\ethm

We already mentioned, back in the opening paragraphs of Section~\ref{S1},
that for a category to be abelian is a property, not a structure. A 
category is abelian if it has finite limits and colimits and these
satisfy certain properties. Applying this to 
$\ct^\heartsuit=\ct^{\leq0}\cap\ct^{\geq0}$,
it follows that there are in $\ct^\heartsuit$
the short exact sequences
that come from its being abelian. When we write $K_i(\ct^\heartsuit)$
we mean the $K$-groups with respect to this natural exact structure.

We should perhaps mention that a sequence $A\la B\la C$, in
the category $\ct^\heartsuit$, is short exact in this
abelian category if and only if there exists, in the 
category $\ct$, a
morphism $C\la\T A$ rendering $A\la B\la C\la\T A$ an exact
triangle.

\section{K-theoretic obstructions to bounded \tstr s}
\label{S4}

We have now arrived at the focus of this survey: we will discuss the
amazing results of Antieau, Gepner and
Heller~\cite{Antieau-Gepner-Heller19}, the many questions they raise,
and the progress since. We begin with the following, which
is a restatement of \cite[Theorems~1.1 and 1.2]{Antieau-Gepner-Heller19}.

\thm{T4.1}
Let $\cm$ be a model category as in \ref{R1.13.3}, and let $\ct=\fgt(\cm)$
be the associated triangulated category. Assume that $\ct$ has a bounded
\tstr\ $\tst\ct$.

Then, with $K_n(\cm)$
as in \ref{R1.13.1}, the following
vanishing results hold:
\be
\item
Unconditionally we have $K_{-1}(\cm)=0$.
\item
If the abelian category $\ct^\heartsuit$ is noetherian, then
$K_n(\cm)=0$ for all $n<0$.
\ee
\ethm

We said at the start that the current survey will highlight the many 
questions that this fascinating theorem raises. We begin with

\plm{P4.3}
The hypotheses in Theorem~\ref{T4.1} are about $\ct=\fgt(\cm)$, and the
conclusions are about $K_i(\cm)$. And there are examples of pairs of 
model categories $\cm$ and $\cn$, with $\fgt(\cm)\cong\fgt(\cn)$ but where
$K_n(\cm)\not\cong K_n(\cn)$.

To elaborate a tiny bit: let $k$ be a field of characteristic $p>0$,
and let $W_2(k)$ be the ring of length-two Witt vectors over $k$. 
For example:
if $k=\zz/p$ then $W_2(k)=\zz/p^2$. One can form the 
singularity categories
$\cm=\mathrm{Sing}\big(W_2(k)\big)$ and 
$\cn=\mathrm{Sing}\big(k[\e]/\e^2\big)$ 
with the natural model
structures, and it is easy to see that they 
satisfy $\fgt(\cm)\cong\fgt(\cn)$. But for
$k=\zz/p$ Schlichting~\cite{Schlichting02}
computes that $K_n(\cm)$ and $K_n(\cn)$ are 
in general non-isomorphic.

Of course: Schlichting's computations are of positive $K$-groups.
For the pair $\cm,\cn$ above
it is easy to show that, for all $n<0$, we have
$K_n(\cm)=0=\K_n(\cn)$.

This immediately raises the question: to what extent do the negative
$K$-groups $K_n(\cm)$ depend on the enhancement $\cm$ of the triangulated
category $\fgt(\cm)$? In particular: is it possible to find a pair
of model categories
$\cm$ and $\cn$, with $\fgt(\cm)\cong\fgt(\cn)$, with $K_{-1}(\cm)=0$ and
with $K_{-1}(\cn)\neq0$? If such a pair exists then the nonvanishing
of $K_{-1}(\cn)$, coupled with Theorem~\ref{T4.1}(i), implies that
the triangulated category $\fgt(\cn)$ cannot have a bounded
\tstr. But the negative $K$-groups of $\cm$ tell us nothing.
\eplm

\rmk{R4.5}
In Problem~\ref{P4.3} we mentioned Schlichting's lovely example.
The triangulated category 
$\D_{\mathbf{sg}}\big(W_2(k)\big)\cong\D_{\mathbf{sg}}\big(k[\e]/\e^2\big)$
has two enhancements with non-isomorphic $K$-groups.

The philosopher can look at this fact in two ways. One could take the
view that this shows that the passage from $\cm$ to $\fgt(\cm)$ loses far
too much information, and one should always work with enhancements.
But one could also take the opposite view,
that there must be something ridiculous
about a theory that distinguishes between the two different enhancements
of 
$\D_{\mathbf{sg}}\big(W_2(k)\big)\cong\D_{\mathbf{sg}}\big(k[\e]/\e^2\big)$.
After all: the triangulated category in question is really dumb. The
objects are the finite-dimensional vector spaces over the field $k$,
the morphisms are the linear maps, the suspension functor is the identity,
and the triangles are the finite 
direct sums of shifts of the triangle
$0\la k\stackrel\id\la k\la 0$.
\ermk

Since this isn't a treatise on philosophy,
let us get back to Theorem~\ref{T4.1} and its ramifications. The negative
$K$-groups, of an enhancement $\cm$ of the triangulated category $\ct$,
carry information about bounded \tstr s on $\ct$. Surely this needs to be
better understood. And so far we have
only raised the question of the dependence
on the enhancement, which for this survey will be a small issue. 

There are two major directions to explore: one can try to improve the 
{\it K--}theoretic results, and one can try to understand better the 
implications for bounded \tstr s. We devote a section to each of those.

\section{K-theoretic generalizations of 
Theorem~\protect{\ref{T4.1}}}
\label{S5}

Possible generalizations
of Theorem~\ref{T4.1} were formulated
already in 
Antieau, Gepner and Heller~\cite{Antieau-Gepner-Heller19}.
We recall

\rmd{R5.1}
With the numbering as in the original
\cite[Introduction]{Antieau-Gepner-Heller19}, the following
conjectures were proposed:
\begin{description}
\item[Conjecture A]
Let $\ca$ be an essentially small abelian category. Then $K_n(\ca)=0$ for
all $n<0$. 
\item[Conjecture B]
Let $\cm$ be an essentially small  model category. If the category $\fgt(\cm)$
has a bounded \tstr\ then $K_n(\cm)=0$ for all $n<0$.
\item[Conjecture C]
Let $\cm$ be an essentially small model category. 
If the category $\fgt(\cm)$
has a bounded \tstr, then the natural map
\[\xymatrix{
K_n\big(\fgt(\cm)^\heartsuit\big)\ar[r] & K_n(\cm)
}\]
is an isomorphism for all $n\in\zz$.
\end{description}
Note that Conjecture~A was not new, it was simply reiterating
the conjecture made by Schlichting more than a decade earlier, see
Reminder~\ref{R2.7}(iv). 

And about Conjecture~C: for $n\geq0$ this is a theorem, not a conjecture.
One can find variants of the result in \cite{Neeman01} and in 
Barwick~\cite{Barwick15}. But both results are about connective
\kth, and say nothing about negative $K$-groups. Thus all three
conjectures are really about negative \kth.

And, after proposing the conjectures in
\cite[Introduction]{Antieau-Gepner-Heller19}, Antieau, Gepner and Heller
go on to observe that Conjecture~B implies the other two conjectures.
\ermd

\dis{D5.3}
Already in  \cite{Antieau-Gepner-Heller19}, Conjecture~B comes with a 
plausibility argument---in fact better still, it comes with the
outline of a proposed proof. And for the purpose of the discussion
that follows I will leave it to the reader to provide
enhancements, the plausibility argument will be stated in terms 
of triangulated categories.

Any time we take an idempotent-complete triangulated category
$\cs$ and a thick subcategory $\car\subset\cs$, we may form the
Verdier quotient $\ct=\cs/\car$. The category $\ct$ will not in 
general be idempotent-complete, but we may form
the idempotent completion $\ct^+$ as in 
Balmer and Schlichting~\cite{Balmer-Schlichting}. And, after enhancing
the picture, we obtain a long exact sequence which, in an abuse of
notation, we write without mention of the enhancements:
\[
\xymatrix@R-10pt@C+20pt{
 & & \ar[dll]\\
K_{-1}(\car) \ar[r]& K_{-1}(\cs)\ar[r] & K_{-1}(\ct^+)\ar[dll]  \\
K_{-2}(\car) \ar[r]& K_{-2}(\cs)\ar[r] & K_{-2}(\ct^+)\ar[dll]  \\
K_{-3}(\car) \ar[r]& K_{-3}(\cs)\ar[r] & K_{-3}(\ct^+)\ar[dll] \\
   & & 
}\]
By Theorem~\ref{T4.1}(i) we know the vanishing of $K_{-1}(\cu)$ for
any $\cu$ with a bounded \tstr. So the idea is to prove the theorem
by dimension shifting.

Suppose therefore that we are given a triangulated category $\car$
with a bounded \tstr. Then it suffices to find an algorithm
that embeds $\car$, as a thick subcategory, in an
idempotent-complete triangulated category $\cs$, and do it in such a 
way that
\be
\item
With $\ct=\cs/\car$, the idempotent completion $\ct^+$ has a 
bounded \tstr.
\item
The inclusion $\car\la\cs$ induces the zero map $K_i(\car)\la K_i(\cs)$.
\ee
The reason this would suffice is that, from the long exact sequence,
we would deduce that the map $K_i(\ct^+)\la K_{i-1}(\car)$ is
surjective. Hence the vanishing of $K_i(\ct^+)$ would permits us
to deduce the vanishing of $K_{i-1}(\car)$.

This plausibility argument persuaded me.
\edis

Perhaps I should elaborate a little. The first I heard of the
beautiful paper~\cite{Antieau-Gepner-Heller19} was when I happened
to be in Boston in April 2018, and Gepner gave a talk about the
results at the MIT topology seminar. Later in 2018 Gepner moved
to Melbourne, and in Australia he gave a number of talks on the same
topic. I attended at least two of these, one in Canberra and one in
Sydney.

And each time, after the talk, I would tell anyone who cared to
listen that this was a lovely paper, and that Conjecture~B had to
be true and the plausibility argument should be turned into a 
proof. After all: we are given a triangulated category $\car$,
with a bounded \tstr, and out of it we want to cook up another
triangulated category $\ct^+$, also with a bounded \tstr. And in
Example~\ref{E0.3} we learned that there are many known recipes for
producing bounded \tstr s. How difficult could it possibly
be to carry out this program?

In early 2020 Covid hit, and in late March 2020 it arrived in 
Australia. And we all went into lockdown. Now: in lockdown there
is not a lot to do, so I decided to try my hand at this myself.
I set out to prove Conjecture~B.  All it required was 
creativity in constructing
bounded \tstr s, and as we have already said: how hard could this
be?

It turns out not to be hard, it's impossible. Being
creative with \tstr s led to the counterexample
of \cite{Neeman21}, which we will
now discuss.

\exm{E5.5}
We need to create a triangulated category with a bounded \tstr.
There are already many known recipes, but we want one lending itself
to computations in 
\kth. Here is one that works.

Let $\ce$ be any idempotent-complete exact category,
and let $\car=\mathrm{Ac}^b(\ce)$ be the category whose objects
are the bounded, acyclic complexes in $\ce$. Recall the definition
in the opening paragraph of \cite[Section~1]{Neeman90}:
a complex
\[\xymatrix{
\cdots\ar[r] & E^{-2} \ar[r] & E^{-1} \ar[r] & E^{0} \ar[r] & 
E^{1} \ar[r] & E^{2} \ar[r] &\cdots
}\]
is declared acyclic if each differential $E^i\la E^{i+1}$ can be
factored as $E^i\la I^i\la E^{i+1}$, with the resulting
$I^{i-1}\la E^i\la I^i$ all admissible exact sequences.
\footnote{
\emph{Caution:} assume $\ce$ is
embedded as an exact subcategory of some abelian category $\ca$.
As in Reminder~\ref{R1.1} this means that a sequence
$E'\la E\la E''$ in $\ce$
is admissible if and only if it is short exact in $\ca$. It
\emph{does not} follow that a cochain complex in $\ce$ of
length $>3$, which
is acyclic in $\ca$, is also acyclic in $\ce$ as
defined in~\cite[Section~1]{Neeman90}. The acyclicity in
the ambient $\ca$
in general depends on the 
embedding.}

This defines the objects in $\car=\mathrm{Ac}^b(\ce)$. The
morphisms in $\car$ are the homotopy equivalence classes of
cochain maps.

And now for the bounded \tstr. In \cite[Lemma~2.1]{Neeman21}
the reader can find a proof that the following is a \tstr\ on
$\mathrm{Ac}^b(\ce)$
\begin{eqnarray*}
\mathrm{Ac}^b(\ce)^{\leq0}&=&\{E^*\in\mathrm{Ac}^b(\ce)\mid
                   E^i=0\text{ for all }i>0\}\ , \\
\mathrm{Ac}^b(\ce)^{\geq0}&=&\{E^*\in\mathrm{Ac}^b(\ce)\mid
                   E^i=0\text{ for all }i<-2\}\ .
\end{eqnarray*}
This \tstr\ is obviously bounded. The problem becomes to compute
the \kth\ of the obvious enhancement.
\eexm

\dis{D5.7}
The strategy of the proof is simple enough. The category 
$\car=\mathrm{Ac}^b(\ce)$ naturally embeds into the category 
$\cs=\K^b(\cs)$, the homotopy category of bounded cochain 
complexes in $\ce$. It is a thick subcategory by 
\cite[Remark~1.10]{Neeman90}, and $\ct=\cs/\car$ is the usual 
derived category
$\D^b(\ce)$. With the standard choice of enhancements
this leads, as in Discussion~\ref{D5.3}, to a long exact
sequence in \kth. And in this case it simplifies to
\[
\xymatrix@R-10pt@C-10pt{
0\ar@{=}[r] &K_{-1}\big(\mathrm{Ac}^b(\ce)\big) \ar[rr]&& K_{-1}(\ce^\oplus)\ar[rr] && K_{-1}(\ce)\ar[dllll] &\, \\
&K_{-2}\big(\mathrm{Ac}^b(\ce)\big) \ar[rr]&& K_{-2}(\ce^\oplus)\ar[rr] && K_{-2}(\ce)\ar[dllll]  \\
& K_{-3}\big(\mathrm{Ac}^b(\ce)\big) \ar[rr]&& K_{-3}(\ce^\oplus)\ar[rr] && K_{-3}(\ce)\ar[dllll] \\
   & & 
}\]
The vanishing is by Theorem~\ref{T4.1}(i).
And the notation, as in Remark~\ref{R1.3},
is that by $\ce^\oplus$ we mean the
category $\ce$ with the split exact structure. Assuming that
Conjecture~B is true and combining with Example~\ref{E5.5},
we have the vanishing of $K_{n}\big(\mathrm{Ac}^b(\ce)\big)$
for all negative $n$. And the long exact sequence tells us that
this would force the natural maps $K_n(\ce^\oplus)\la K_n(\ce)$
to be isomorphisms for all $n<0$. 

We will be sketching one approach that leads to counterexamples.
\edis

\dis{D5.9}
If we are going to kill a conjecture, then why not kill two with one stone?

With this philosophy in mind, let $\ca$ be the abelian category
$\ca=\mathrm{Ac}^b(\ce)^\heartsuit$.
If we could guarantee that the natural map 
$K_i(\ca)\la K_i\big(\mathrm{Ac}^b(\ce)\big)$ is an
isomorphism, then an $\ce$ providing a counterexample to Conjecture~B
would automatically also kill Conjecture~A. In the light of \ref{R1.13.3},
it suffices to prove that the natural map $\D^b(\ca)\la\mathrm{Ac}^b(\ce)$
is an equivalence of triangulated categories.

This isn't true for every $\ce$. The matter is studied in
\cite[Section~2]{Neeman21}, and \cite[Proposition~2.4]{Neeman21}
shows that
the natural map $\D^b(\ca)\la\mathrm{Ac}^b(\ce)$
is an equivalence
if and only if the category $\ce$ is 
hereditary. This means: if and only if, for every pair of objects
$E,E'\in\ce$, we have $\Ext^n(E,E')=0$ for all $n>1$.

If $X$ is an algebraic curve, and $\ce=\vect X$ is the exact
category of vector bundles over $X$, then the category $\ce$ is hereditary.
\edis

\smr{S5.11}
With Example~\ref{E5.5} and Discussions~\ref{D5.7} and
\ref{D5.9} in mind, the article \cite{Neeman21} restricts attention
to the exact category $\ce=\vect X$ with $X$ a projective curve over
a field $k$.
The article then proves 
\be
\item
If the singularities of $X$ are no worse than simple nodes, then
$K_{-1}(\ce^\oplus)=0$.
\setcounter{enumiv}{\value{enumi}}
\ee
And this is combined with an appeal to the classical literature, after 
all 
\be
\setcounter{enumi}{\value{enumiv}}
\item
There are known examples of nodal curves $X$ such that, with $\ce=\vect X$,
we have $K_{-1}(\ce)\neq0$.
\ee
The conjunction of (i) and (ii) with the earlier discussion 
kills Conjectures~A and B. 
\esmr

\rmk{R5.13}
Why is it that, after the paper is published and can no longer be changed,
we always wish we could edit it? Before I explain the part I would change
let me make two points.
\be
\item
There is an inevitability about it. Whatever argument we come
up with, and no matter how clever we find it at the time, someone
will eventually discover an approach that is
shorter, simpler, more elegant and 
more general. It could be us, or it might be someone else. It can
happen quickly or take years. One thing is guaranteed: the
passage of time will bring improvements.
\setcounter{enumiv}{\value{enumi}}
\ee
As the reader will discover, the part I would rewrite if I could is
contained in Summary~\ref{S5.11}(i).
In my defense let me say that I found the statement difficult to believe.
In the coming paragraphs I will try to explain my skepticism.

As explained in Discussion~\ref{D5.7}: the search for
a counterexample to Conjecture~B 
was reduced to looking for exact categories
$\ce$, with $K_n(\ce^\oplus)\not\cong K_n(\ce)$ for some $n<0$.
To kill Conjectures~A and B together, Discussion~\ref{D5.9}
tells us that a hereditary $\ce$
is preferred.
And vector bundles over projective curves were 
a reasonable first place to look, as we will now explain.

We have mentioned that, on any curve $X$, the 
exact category $\vect X$
is hereditary.
An irreducible curve has two choices: it is either affine or projective.
For affine curves $X=\spec R$, we have that $\vect X$ is equivalent
to the category $\proj R$ of finitely generated, projective $R$-modules.
And in the category $\ce=\proj R$ every short exact sequence splits,
that is $\ce^\oplus=\ce$. For us this eliminates affine curves, on 
affine curves the map $K_n(\ce^\oplus)\la K_n(\ce)$ 
is guaranteed to
be an isomorphism. 

Projective curves $X$ have non-split short exact sequences of vector
bundles, hence studying $\ce=\vect X$ is natural enough. What surprised me
wasn't that $K_{-1}(\ce^\oplus)\not\cong K_{-1}(\ce)$, but that
$K_{-1}(\ce^\oplus)=0$.
And the cause for my surprise was that, for affine curves $X$ and with
$\ce=\vect X$, we often have 
$K_{-1}(\ce^\oplus)=\K_{-1}(\ce)\neq0$. What made projective curves
so totally unlike affine curves? I was expecting some difference,
but not such a drastic, radical dichotomy.

Because the result seemed suspect to me, I checked it carefully.
My second excuse for my stupidity is that
\be
\setcounter{enumi}{\value{enumiv}}
\item
In the process of carefully proving that a result is true, we can
lose sight of the bigger picture. We fail to ask ourselves for the
underlying reason. The all-important question should be:
\emph{Why} is this surprising result true?
\ee
The time has come to stop making excuses, and tell the reader how I would
rewrite the article if I could.
\ermk

If the Almighty Lord granted me a time machine, and permitted 
me to rewrite~\cite{Neeman21} 
and have the improved version published in place of
the one currently in print, then instead of
Summary~\ref{S5.11}(i) the result would be:

\thm{T5.15}
Let $\ca$ be an idempotent-complete additive category. Assume that,
for every object $A\in\ca$, the ring $\Hom(A,A)$ is artinian.
Then $K_n(\ca^\oplus)=0$ for all $n<0$.
\ethm

\prf
It is known that, if $R$ is any artinian ring, then $K_n(R)=0$ for
all $n<0$; for a reference the reader can see 
Weibel~\cite[Theorem~2.3]{Weibel89A}. Rephrasing this in terms of
categories: the category $\cb=\proj R$, of finitely generated,
projective $R$-modules, has vanishing negative \kth.

Now choose an object $A\in\ca$, Applying the above to the
ring $R=\Hom(A,A)$ gives that the full,
additive subcategory $\add(A)\subset\ca$,
of all direct summands of
finite direct sums $A^{\oplus m}$ of $A$ with itself, is such that
$K_n\big(\add(A)^\oplus\big)=0$ for all $n<0$. But $\ca$ is the filtered
union of the full
subcategories $\add(A)$, making 
$K_n(\ca^\oplus)$ the direct limit of
the $K_n\big(\add(A)^\oplus\big)$. And
for $n<0$ this yields $K_n(\ca^\oplus)=0$.
\eprf

\plm{P5.17}
The emphasis in this survey is on open problems, and Conjecture~C is 
still open.
Theorem~\ref{T5.15} might be relevant, it produces many examples
of model categories with bounded \tstr s, and whose negative \kth\ 
can be
computed and is nonzero.

Let $X$ be any projective scheme and put $\ce=\vect X$. By 
Theorem~\ref{T5.15} we know that $K_n(\ce^\oplus)=0$ for all
$n<0$, and the exact sequence of Discussion~\ref{D5.7}
tells us that the map $K_n(\ce)\la K_{n-1}\big(\text{Ac}^b(\ce)\big)$
is an isomorphism for all $n<0$. And there are many examples out
there of projective schemes for which $K_n(\ce)$ has been computed
when $n<0$, and is nonzero.

Of course: I have no idea how to go about computing 
$K_{n}\big(\text{Ac}^b(\ce)^\heartsuit\big)$, except in the case 
where $\dim(X)=1$ and 
$K_{n}\big(\text{Ac}^b(\ce)^\heartsuit\big)$ 
agrees with $K_{n}\big(\text{Ac}^b(\ce)\big)$
by Discussion~\ref{D5.9}.
But since this provides such a wealth of candidates, one should
at least do a reality check on the conjecture. If $X$ and $Y$ are two
projective schemes, we can put $\ce=\vect X$ and $\cf=\vect Y$.
The above gives a procedure for computing
 $K_{n}\big(\text{Ac}^b(\ce)\big)$ and  
$K_{n}\big(\text{Ac}^b(\cf)\big)$. If Conjecture~C is true, then 
these depend only on the heart of the bounded \tstr s. Thus an 
equivalence
\[
\text{Ac}^b(\ce)^\heartsuit\cong \text{Ac}^b(\cf)^\heartsuit
\]
would imply that 
$K_{n}\big(\text{Ac}^b(\ce)\big)\cong K_{n}\big(\text{Ac}^b(\cf)\big)$
for all $n$. And for $n<0$ the discussion shows that this
means $K_n(\ce)\cong K_n(\cf)$. Or in more classical notation:
an equivalence of the hearts of these \tstr s would force
$K_n(X)\cong K_n(Y)$ for all $n<0$.

This should be checked.
\eplm

\section{Categories that can be proved to have no
bounded \tstr}
\label{S6}

In Section~\ref{S5} the focus was on conjectures aimed at strengthening
the vanishing statements in negative \kth. In the current section
we will ask if the results about bounded \tstr s can be sharpened.

The title of Antieau, Gepner and 
Heller's article~\cite{Antieau-Gepner-Heller19}
gives the clue: Theorem~\ref{T4.1}, which is a restatement
of \cite[Theorems~1.1 and 1.2]{Antieau-Gepner-Heller19},
provides obstructions to the existence of bounded \tstr s. 
As we mentioned already in Problem~\ref{P0.1}: there are triangulated
categories of motives which are conjectured to have bounded \tstr s
and, armed with Theorem~\ref{T4.1}, we can test these conjectures
with a reality check. If their negative \kth\ fails to vanish then
either they have no bounded \tstr s, or these \tstr s must have
non-noetherian hearts.

Over the decades there has been much work computing the negative \kth\ 
of schemes, and there are many examples out there of schemes with
nonvanishing negative \kth. This immediately leads to 
\cite[Corollary~1.4]{Antieau-Gepner-Heller19}. It says:

\cor{C6.1}
Let $X$ be a noetherian scheme. Then the following holds:
\be
\item
If $K_{-1}(X)\neq0$ then $\dperf X$ has no bounded \tstr.
\item
If $K_n(X)\neq0$ for $n<-1$, then a bounded \tstr\ on $\dperf X$
cannot have a noetherian heart.
\ee
\ecor

\nin
Being courageous, Antieau, Gepner and 
Heller go on to 
formulate \cite[Conjecture~1.5]{Antieau-Gepner-Heller19}. It
predicts

\cnj{C6.3}
Let $X$ be a finite-dimensional, noetherian scheme. Then the category
$\dperf X$ has a bounded \tstr\ if and only if $X$ is regular.
\ecnj

\dis{D6.5}
\cite[Conjecture~B]{Antieau-Gepner-Heller19}
was eminently 
plausible and came 
with a projected outline of a possible proof. We
sketched it in Discussion~\ref{D5.3}.

The evidence for 
\cite[Conjecture~1.5]{Antieau-Gepner-Heller19} was thin  by
comparison. There are many singular noetherian schemes whose
negative \kth\ vanishes. Recall that all zero-dimensional,
noetherian
schemes have vanishing negative \kth---we met this in the proof
of Theorem~\ref{T5.15}, where we appealed to a non-commutative 
generalization. And it is easy to construct
zero-dimensional, noetherian schemes with dreadful singularities.
In this light, it was a bold leap of faith for
Antieau, Gepner and 
Heller to pass from the meager evidence, provided by
Corollary~\ref{C6.1}, to the daring Conjecture~\ref{C6.3}.
\edis

Boldness can pay off, this conjecture turns out to be very true---meaning so true that
a vast generalization is now a theorem.
But let us proceed chronologically: the first progress
on the conjecture came in Harry Smith~\cite{SmithHarry19}: Smith proved
the conjecture for affine $X$. For affine $X$ there is a
classification of the compactly generated \tstr s
on $\Dqc(X)$ due to Alonso, Jerem{\'{\i}}as and
Saor{\'{\i}}n~\cite{Alonso-Jeremias-Saorin10},
and an analysis of which of these compactly generated
\tstr s restrict to \tstr s on $\dcoh(X)$. Using this machinery,
Smith was able to prove Conjecture~\ref{C6.3} for affine $X$.

In fact Smith proved a stronger statement: if $X$ is singular,
affine,
irreducible, noetherian and
finite-dimensional then the only \tstr s on $\dperf X$ are the
trivial ones. This means that, under the hypotheses above,
either
${\dperf X}^{\leq0}=0$ or
${\dperf X}^{\geq0}=0$.
Note that this \emph{does not} generalize beyond
the affine case. There are known examples of singular varieties
for which $\dperf X$ has nontrivial semiorthogonal
decompositions---which are of course examples
of non-trivial \tstr s. And when one of the
components is admissible, \ref{E0.3.1} allows us to
glue \tstr s on the components, producing even
more non-trivial \tstr s.
There is a rich and growing literature on the
possible, nontrivial semiorthogonal
decompositions of $\dperf X$. For smooth $X$ this
literature is immense and
goes back decades, it originates with
Be{\u\i}linson~\cite{Beilinson83} and the semiorthogonal
decomposition of $\dperf{\pp^n}$; in this case
all of the components are admissible and we deduce
an abundance of bounded \tstr s.
For singular $X$, the study of semiorthogonal decompositions
of $\dperf X$ is a much more recent development, but
is an active field: we refer the
reader to
Kalck, Pavic and Shinder~\cite{Kalck-Pavic-Shinder21}
and to
Karmazyn, Kuznetsov and
Shinder~\cite{Karmazyn-Kuznetsov-Shinder22}
for a tiny sample.

Semiorthogonal decompositions of
$\dperf X$ can give exotic \tstr s,
both bounded and unbounded. In view of
this any classification theorem,
of the possible \tstr s
on $\dcoh(X)$ or $\dperf X$ for $X$ non-affine,
would have to be much more delicate than the
analysis in
Alonso, Jerem{\'{\i}}as and
Saor{\'{\i}}n~\cite{Alonso-Jeremias-Saorin10}.
Harry Smith's approach, which proved
Conjecture~\ref{C6.3} in the affine case
using the results of \cite{Alonso-Jeremias-Saorin10},
does not generalize in any straightforward way.
Nevertheless the conjecture is
true unconditionally---not only is it
true as originally stated, 
but so is a major generalization. The following
is now a theorem,
see \cite[Theorem~0.1]{Neeman22A}.

\thm{T6.7}
Let $X$ be a finite-dimensional, noetherian scheme, and let $Z\subset X$
be a closed subset. Let $\dperfs ZX$ be the derived category, whose
objects are the perfect
complexes on $X$ which restrict to acyclic complexes on
the open set $X-Z$.

The category $\dperfs ZX$ has a bounded \tstr\ if and only if $Z$ 
is contained in the regular locus of $X$.  
\ethm

\rmk{R6.9}
The proof of Theorem~\ref{T6.7} is not by computing an
obstruction. It uses the metric techniques introduced by the author
in the theory of approximable triangulated categories. For a survey
of that theory the reader is referred for example to~\cite{Neeman22}.
The current survey is not the right place to expand much on this 
cryptic comment---in this remark we will confine ourselves to
a bare minimum.

Let
$\Dqcs Z(X)$ be the derived category, whose objects are
complexes of $\co_X^{}$-modules with quasicoherent cohomology
supported on $Z$. This category has a single compact generator,
hence a preferred equivalence class of \tstr s. Take any
\tstr\ in this preferred equivalence class and view it as
a metric on $\Dqcs Z(X)$, in the sense explained in
\cite{Neeman22}. This restricts to a metric on the subspace
$\dperfs ZX\subset\Dqcs Z(X)$, and up to equivalence these
metrics are independent of any choice. We have a preferred
equivalence
class of metrics on $\dperfs ZX$.

The key to proving Theorem~\ref{T6.7} turns out to be showing
that any bounded \tstr\ on $\dperfs ZX$ yields a metric in this
preferred equivalence class. And, if the reader will permit us to
use yet more of the machinery exposed in \cite{Neeman22}, there is
a process which passes from one triangulated category with a metric to
another, it constructs out of $\cs$ a category $\fs(\cs)$. This
can of course be applied to $\cs=\dperfs ZX$ with any of the
equivalent  metrics above. And  
proof then analyses the category $\fs\big(\dperfs ZX\big)$, and
shows that it agrees with $\dperfs ZX$.
\ermk

\plm{P6.9}
Theorem~\ref{T6.7} is proved by methods which don't fit into
this survey. But, now that we know the theorem to be true, it does
raise a natural question.

Clearly negative \kth\ is not the only obstruction to the 
existence of bounded \tstr s. Given how important bounded \tstr s
have turned out to be, across wide swaths of
mathematics, what is the right obstruction? Can one
find an easily computable obstruction to the existence of
bounded
\tstrs, which doesn't vanish on categories such as
$\dperfs ZX$?  

And once again we remind the reader of Problem~\ref{P0.1}.
By now there
exist a plethora of derived categories of motives, all conjectured
to have bounded \tstr s. An obstruction that is easily computable 
on these candidates would be a treasure, it would allow us to
weed out the ones with no chance of working. 
\eplm

\providecommand{\href}[2]{#2}

\end{document}